\documentclass[12pt,oneside, italian]{article}
\usepackage{latexsym}
\usepackage{color}
\usepackage{amssymb,amsthm,amsmath}
\usepackage{graphicx}

\newtheorem{theorem}{Theorem}

\newtheorem{conjecture}{Conjecture}
\begin{document}
\title{The maximum and the minimum size\\ of complete $(n,3)$-arcs in $PG(2,16)$}
\date{}
\author{D. Bartoli, S. Marcugini and F. Pambianco \\
%EndAName
Dipartimento di Matematica e Informatica,\\
Universit\`{a} degli Studi di Perugia, \\
Via Vanvitelli 1, 06123 Perugia Italy \\
e-mail: \{daniele.bartoli, gino, fernanda\}@dmi.unipg.it}
\maketitle
\setlength{\baselineskip}{15 pt}
\begin{abstract}
\noindent In this work we solve the packing problem for complete $(n,3)$-arcs in $PG(2,16)$, determining that the maximum size is $28$  and the minimum size is $15$. We also performed a partial classification of the extremal size of complete $(n,3)$-arcs in $PG(2,16)$.
\end{abstract}
%{\em Key words: }Covering codes, Projective geometry, Saturating sets

\section{Introduction}
In the projective plane $PG(2,q)$ over the finite field $GF(q)$ an $(n,r)$-arc is a set of
$n$ points such that no $r + 1$ points are collinear and some $r$ points are collinear. An $(n,r)$-arc is called \emph{complete} if it is not contained in a $(n+1,r)$-arc of the same projective plane. An $(n,2)$-arc is called $n$-arc. For a more detailed introduction to $(n,r)$-arcs and in particular $(n,3)$-arcs see \cite{H1988}, \cite{HS2000}, \cite{3ArchiMPG2_11},  \cite{3ArchiMPG2_13}. The largest size of $(n,r)$-arcs of $PG(2,q)$ is indicated by $m_{r}(2,q)$. In particular $m_{3}(2,q)\leq 2q+1$ for $q\geq 4$ (see \cite{Thas1975}). In \cite{HS2000} bounds for $m_{r}(2,q)$ and the relationship between the theory of complete $(n,r)$-arcs, coding theory and mathematical statistics are given.\\
Arcs and and $(n,3)$-arcs in $PG(2,q)$ correspond to respectively MDS and NMDS codes of dimension $3$. These types of linear codes are the best in term of minimum distance, among the linear codes with the same length and dimension. In general $(n,k)$-arcs in $PG(2,q)$ correspond to linear codes with Singleton defect equal to $k-2$.\\

\section{Results}
In this work we establish the maximum and the minimum size of complete $(n,3)$-arcs in $PG(2,16)$. To do this we performed a computer based search using some ideas similar to those presented in \cite{3ArchiPG2_7}, \cite{3ArchiMPG2_11}, \cite{3ArchiMPG2_13} and \cite{Cook2011}.
\begin{theorem}
The maximum size of complete $(n,3)$-arcs in $PG(2,16)$ is $28$.
\end{theorem}
\proof
We performed an exhaustive search of $(n,3)$-arcs in $PG(2,16)$ of size greater than $28$ and we found no examples. Moreover we have an example of complete $(28,3)$-arc (see \cite{BFMP2006}), obtained as union of orbits of some subgroup of $P\Gamma L(3,16)$; the size is equal to the one given in \cite{BKW2005}.
\endproof
The classification of complete $(28,3)$-arcs is in progress.

\begin{theorem}
The minimum size of complete $(n,3)$-arcs in $PG(2,16)$ is $15$.
\end{theorem}
\proof We performed an exhaustive search of $(n,3)$-arcs in $PG(2,16)$ of size less than
$15$ and we found no examples. We also proved that a complete $(15,3)$-arc contains a
$(k,2)$-arc, with $8\leq k \leq 9$. As result of the search for complete $(15,3)$-arcs we
get only the example presented in Table \ref{15arco} containing a $(9,2)$-arc.
\endproof

\begin{conjecture}
There exists a unique complete $(15,3)$-arc in $PG(2,16)$.
\end{conjecture}

We denote{\bf \ }$GF(16)=\{0,1=\alpha ^{0},2=\alpha ^{1},\ldots ,15=\alpha ^{14}\}$ where $%
\alpha $ is a primitive element such that $\alpha^4 + \alpha^3 + 1=0$. The columns $\ell_{i}$ indicate the number of $i$-secant of the $(n,3)$-arc and $G$ indicates the description of the stabilizer in $P\Gamma L(3,16)$ (see \cite{librogruppi}).

\begin{center}
\begin{table}
\caption{complete $(15,3)$-arc }\label{15arco}
\begin{center}
\begin{tabular}{|c|c|c|c|c|c|}
\hline
Points &$\ell_{0}$&$\ell_{1}$&$\ell_{2}$&$\ell_{3}$&G\\
\hline
\hline
\tabcolsep=0.85 mm
\begin{tabular}{ccccccccccccccc}
 1& 0& 0& 1& 1& 1& 1& 1& 1& 1& 1& 1& 1& 1& 1\\
 0& 1& 0& 1& 0& 1& 2& 2& 4& 9& 9& 11& 11& 13& 13\\
 0& 0& 1& 1& 11& 8& 5& 10& 10& 2& 8& 2& 11& 1& 12\\
\end{tabular}&
92 &
138 &
12 &
31 &
$\mathcal{S}_3$\\
\hline
\end{tabular}
\end{center}
\end{table}
\end{center}

\end{document}